
\magnification=1200
\overfullrule=0mm
\nopagenumbers
\vskip 1cm

\centerline {\bf On an Algebraic System of Equations}
\bigskip
\centerline {\bf H. Hakopian, M. Tonoyan}
\bigskip

We study the following algebraic system of equations
$$
P_\alpha (x_1,\ldots ,x_n):=P_\alpha (x):=
x^\alpha -  \sum _{\beta \in I} a_{\alpha , \beta}x^\beta = 0,
\quad \alpha \in J, \eqno (1)
$$
where
$I = \{\alpha \in Z_+^n:|\alpha|\le m\},\quad
J:= J_{m+1} = \{\alpha \in Z_+^n:|\alpha| = m+1\}.$

\noindent We consider some $\# I \times \# I$ - matrices:
$A_1,\ldots ,A_n$.
The structure of these matrices is as follows: $k$ rows ($k = \# J_m$) are
the coefficients of $k$ polynomials from
$\{P_\alpha: \alpha \in J \}$
and each of the remaining rows consists of $0$'s and one $1$.

\noindent We prove that the number of distinct solutions to (1) equals to
maximal possible $\ \# I\ $
if and only if the above matrices are commuting
and each of them is semisimple (i.e. they have complete system of
eigenvectors).

This gives a characterization of correct systems of points for
multivariate Lagrange interpolation.
This also reduces the solving of (1) to solving of $n$ univariate algebraic
equations: to finding the eigenvalues of the above matrices.
In the case when the eigenvalues of some matrix are distinct,
the solving of (1) is reducing to solving of only one univariate
equation: to finding of these distinct eigenvalues.

The above results remain valid also for more general systems (1):
with $I$ being any lower set.
                                  
\bye